\documentclass[11pt]{amsart}

\usepackage{amsmath,amssymb,amsthm,amsrefs,a4wide}
\usepackage{latexsym}
\usepackage{indentfirst}
\usepackage{graphicx}
\usepackage{placeins}
\usepackage{booktabs}
\usepackage{algorithm}
\usepackage{algorithmic}
\usepackage{multirow}
\usepackage{color}

\numberwithin{figure}{section}

\theoremstyle{plain}

\theoremstyle{definition}

\newtheorem{example}{Example}[section]

\theoremstyle{remark}

\newcommand{\diag}{\mathrm{diag}}

\newcommand{\eps}{\epsilon}

\newcommand{\bbm}{\begin{bmatrix}}
\newcommand{\ebm}{\end{bmatrix}}
\newcommand{\R}{\mathbb{R}}
\newcommand{\C}{\mathbb{C}}

\newcommand{\T}{\mathsf{T}}

\newcommand{\nx}{{n}}
\newcommand{\nz}{{n_z}}
\newcommand{\nc}{{n_c}}
\newcommand{\nl}{{n_l}}

\newcommand{\Bh}{\hat{B}}

\newcommand{\bb}{\mathbf{b}}
\newcommand{\bu}{\mathbf{u}}

\newcommand{\bbh}{\hat{\bb}}

\newcommand{\muAa}{{\mu_{A_1}}}
\newcommand{\muAb}{{\mu_{A_2}}}
\newcommand{\muAn}{{\mu_{A_n}}}
\newcommand{\muC}{{\mu_C}}

\begin{document}

\title[Blind Free Deconvolution over one-parameter sparse families via eigenmatrix]{Blind Free Deconvolution over one-parameter sparse families via eigenmatrix}

\author[]{Lexing Ying} \address[Lexing Ying]{Department of Mathematics, Stanford University,
  Stanford, CA 94305} \email{lexing@stanford.edu}

\thanks{This work is partially supported by NSF grant DMS-2208163. The author thanks Xiucai Ding for helpful discussions.}



\begin{abstract}
  This note considers the blind free deconvolution problems of sparse spectral measures from one-parameter families. These problems pose significant challenges since they involve nonlinear sparse recovery. The main technical tool is the eigenmatrix method for solving unstructured sparse recovery problems. The key idea is to turn the nonlinear inverse problem into a linear inverse problem by leveraging the R-transform for free addition and the S-transform for free product. The resulting linear problem is solved with the eigenmatrix method tailored to the domain of the parametric family. Numerical results are provided for both the additive and multiplicative free deconvolutions.
\end{abstract}

\maketitle

\section{Introduction}\label{sec:intro}

This note considers two blind free deconvolution problems over sparse spectral measures from one-parameter families. In the additive setting, $A_1$ and $A_2$ are both $N\times N$ symmetric matrices with spectral measures $\muAa$ and $\muAb$. Given the spectral measure $\muC$ of $C = A_1 + Q A_2 Q^\T$, where $Q$ is sampled from the Haar measure over orthogonal matrices, the task is to recover $\muAa$ and $\muAb$. This problem can also be extended to the free sum of more than two terms: given $\muC$ of $C = A_1 + Q_2 A_2 Q_2^\T + \ldots + Q_n A_n Q_n^\T$, the task is to recover the spectral measures $\mu_{A_1}, \ldots, \mu_{A_n}$.

In the multiplicative setting, $A_1$ and $A_2$ are both $N\times N$ symmetric positive definite matrices with sparse spectral measures. Given $\muC$ of $C = \sqrt{A_1} Q A_2 Q^\T \sqrt{A_1}$, where $Q$ is again from the Haar measure over orthogonal matrices, the task is to recover the spectral measures $\muAa$ and $\muAb$. This problem can also be extended to the free product of more than two terms: given $\muC$ of $C =\sqrt{A_1} Q_2 \sqrt{A_2} \ldots Q_n A_n Q_n^\T \ldots \sqrt{A_2} Q_2^\T \sqrt{A_1}$, the task is to recover the spectral measures $\mu_{A_1}, \ldots, \mu_{A_n}$.

These problems are referred to as blind free deconvolution since neither $\muAa$ nor $\muAb$ is known. Blind free deconvolutions for general spectral measures are ill-posed. For example, the sum of two semicircle laws is a semicircle law with the sum of their widths, and therefore, it is not possible to recover the individual inputs. In this note, we restrict to sparse spectral measures. Even with the sparsity constraint, the problem remains difficult since the composition of the spectral measures of $A_1$ and $A_2$ into that of $C$ is nonlinear, thereby rendering these two problems nonlinear inverse problems.

\subsection{Related work.} There is a significant body of work on spectral estimation, including the linear shrinkage \cite{ledoit2004well}, the optimization approach using the Marcenko-Pastur equation \cite{el2008spectrum}, the nonlinear shrinkage methods \cite{ledoit2011eigenvectors,ledoit2012nonlinear,ledoit2015spectrum,ledoit2020analytical,gavish2014optimal,bun2016rotational}, the moment-based method \cite{kong2017spectrum}, and the subordination method \cite{arizmendi2020subordination}. In most of the previous work, the estimation task is not of the blind deconvolution type, as one of the convolution components (such as the noise level or the dimension-to-sample size ratio) is assumed to be known. Relatively few works have focused on blind free deconvolution. Recently, the work in \cite{ying2025sparse} considers the spectral estimation with an unknown noise level, which can be viewed as a relatively simpler version of blind deconvolution.

\subsection{Contributions.} This note takes some preliminary steps in investigating these two blind deconvolution problems. In addition to the sparsity of the spectral measures, we further restrict to a simple setting such that $\mu_{A_k}$ is from certain one-parameter families.

The primary technical tool is the recently proposed eigenmatrix method \cite{ying2024eigenmatrix,ying2025multidimensional} for solving general, yet unstructured, sparse recovery problems. This level of generality will play a key role in addressing the problems here.

There are two key ideas. First, the nonlinear inverse problems can be turned into linear ones by leveraging the R-transform for free addition and the S-transform for free product \cite{voiculescu1986addition,voiculescu1987multiplication,mingo2017free,potters2020first}. Once this is done, we extend the eigenmatrix method to the parametric case to solve the corresponding linear sparse recovery problems.

The rest of the note is organized as follows. Section \ref{sec:em} reviews the eigenmatrix method.
Section \ref{sec:cls} presents the main ideas using classical deconvolution as a warm-up.
Section \ref{sec:add} describes the additive case.
Section \ref{sec:mul} discusses the multiplicative case.
Section \ref{sec:disc} concludes with a discussion for future work.

\section{Eigenmatrix} \label{sec:em}

This section provides a general introduction to the eigenmatrix method for unstructured sparse recovery problems. Let $X$ be the parameter space and $G(z,x)$ be a kernel function defined for $x\in X$ at sample $z$. We assume that $G(z,x)$ is analytic in $x$. Throughout this note, $X$ is an interval of $\R$, $z\in \C$, but $G(z,x)$ could be quite general. Suppose that
\begin{equation}
  f(x) = \sum_{k=1}^{\nx} w_k \delta(x-x_k)
  \label{eq:f}
\end{equation}
is the unknown sparse signal, where $\{x_k\}_{1\le k \le \nx}$ are the spike locations, $\{w_k\}_{1\le k \le \nx}$ are the spike weights, and $\delta(\cdot)$ is the Dirac delta function at the origin. The observed data is defined via the function
\begin{equation}
  u(z) :=\int_X G(z,x) f(x) dx = \sum_{k=1}^{\nx} G(z,x_k) w_k.
  \label{eq:u}
\end{equation}

Let $\{z_j\}$ be a set of $\nz$ unstructured samples and $u_j \approx u(z_j)$ be the noisy observations of $u(z)$ at the chosen samples. The task is to recover the spikes $\{x_k\}$ and weights $\{w_k\}$ from $\{u_j\}$.

Define for each $x$ the column vector
\[
\bb_x:= 
\bbm
G(z_j,x)
\ebm_{1\le j \le \nz}
\]
in $\C^{n_z}$. As a function of $x$, $\bb_x$ is also analytic in $x$.

The first step is to construct a matrix $M \in \C^{n_z\times n_z}$ such that $M \bb_{x} \approx x \bb_{x}$ for $x\in X$, i.e., $M$ is the matrix with $(x,\bb_{x})$ as approximate eigenpairs for $x\in X$.  Numerically, it is more robust to use the normalized vector $\bbh_{x} = \bb_{x}/\|\bb_{x}\|$ since the norm of $\bb_{x}$ can vary significantly depending on $x$. The condition then becomes
\[
M \bbh_{x} \approx x \bbh_{x}, \quad x\in X.
\]
To build $M$, we pick a Chebyshev grid $\{c_t\}_{1\le t \le \nc}$ of size $\nc$  \cite{trefethen2019approximation} of the interval $X$. Here $\nc$ is sufficiently large yet the vectors $\{\bbh_{c_t}\}$ are numerically linearly independent. This condition is then enforced on this Chebyshev grid, i.e.,
\[
M \bbh_{c_t} \approx c_t \bbh_{c_t}.
\]
By defining the $\nz\times \nc$ matrix
\begin{equation}
  \Bh = 
  \bbm
  \bbh_{c_1} & \ldots & \bbh_{c_\nc}
  \ebm
  \label{eq:Bh}
\end{equation}
and the $\nc\times \nc$ diagonal matrix $\Lambda = \diag(c_t)$, the previous condition can be written as
\[
M \Bh \approx \Bh \Lambda.
\]
Under the assumptions that $\bb_{x}$ is analytic (or smooth) in $x$ and the columns of $\Bh$ are numerically linearly independent, we define the eigenmatrix $M$
\begin{equation}
  M := \Bh \Lambda \Bh^+, \label{eq:M}
\end{equation}
where the pseudoinverse $\Bh^+$ is computed by thresholding the singular values of $\Bh$. In practice, the thresholding value is chosen so that a small constant bounds the operator norm of $M$. 

The rest of the method follows, for example, the ESPRIT algorithm \cite{roy1989esprit} (or the Prony method \cite{prony1795essai}). Define the vector
\[
\bu=
\bbm   
u_j
\ebm_{1\le j \le \nz}
\]
in $\C^{n_z}$, where $u_j$ are the noisy observations. Notice that $\bu = \sum_k \bb_{x_k} w_k$. Consider the matrix
\[
T \equiv \bbm
\bu & M\bu & \ldots & M^\nl \bu
\ebm
\]
with $\nl > \nx$, obtained from applying $M$ repetitively. Since $\bu \approx \sum_k \bb_{x_k} w_k$ and $M\bb_{x}\approx x \bb_{x}$,
\[
T = \bbm  \bu & M\bu & \ldots & M^\nl \bu \ebm
\approx
\bbm \bb_{x_1} & \ldots & \bb_{x_\nx} \ebm
\bbm w_1 & &\\& \ddots & \\& & w_\nx\ebm
\bbm
1 &  x_1 & \ldots & (x_1)^\nl \\
\vdots & \vdots & \ddots & \vdots \\
1 &  x_\nx & \ldots & (x_\nx)^\nl
\ebm.
\]
Let $U S V^*$ be the rank-$\nx$ truncated SVD of $T$. The matrix ${V}^*$ satisfies
\[
{V}^* \approx P
\bbm
1 &  x_1 & \ldots & (x_1)^\nl \\
\vdots & \vdots & \ddots & \vdots \\
1 &  x_\nx & \ldots & (x_\nx)^\nl
\ebm,
\]
where $P$ is an unknown non-degenerate $\nx \times \nx$ matrix. Let $Z_L$ and $Z_H$ be the submatrices obtained by excluding the last column and the first column of $V^*$, respectively, i.e.,
\[
Z_L \approx P
\bbm
1      & \ldots & (x_1)^{\nl-1} \\
\vdots & \ddots & \vdots \\
1      & \ldots & (x_\nx)^{\nl-1}
\ebm,
\quad
Z_H \approx P
\bbm
x_1 & \ldots & (x_1)^\nl \\
\vdots & \ddots & \vdots \\
x_\nx & \ldots & (x_\nx)^\nl
\ebm.
\]
By forming $Z_H (Z_L)^+$ and noticing
\[
Z_H (Z_L)^+ \approx
P
\bbm
x_1 &  & \\
& \ddots & \\
& & x_\nx
\ebm
P^{-1},
\]
one obtains the estimates for $\{x_k\}$ by computing the eigenvalues of $Z_H (Z_L)^+$.

With the estimates for $\{x_k\}$ available, the least square solution of 
\[
\min_{w_k} \sum_j \left|\sum_k G(s_j,x_k) w_k - u_j\right|^2 
\]
gives the estimators for $\{w_k\}$.

\section{Classical deconvolution}\label{sec:cls}

This section illustrates the main ideas using the classical deconvolution as an example, which is a new application of the eigenmatrix method. Let the interval $X$ be the parametric domain of the 1-parameter family of distributions $\rho_x(y)$ over $\R$. Each distribution $\alpha_k$ is associated with the parameter $x_k$, i.e., $\alpha_k(y) = \rho_{x_k}(y)$.  Given $\gamma = \alpha_1 * \ldots * \alpha_n$, the task is to recover the parameters $x_1,\ldots,x_n$, and hence $\alpha_1, \ldots, \alpha_n$.

The first step turns this nonlinear deconvolution problem into a linear one. Consider the characteristic function
\[
\hat{\alpha}_k (\xi) = \int e^{-i y\xi } d \alpha_k(y),
\]
for each $\alpha_k$ and similarly for $\gamma$. Then, $\hat{\alpha}_1(\xi) \cdot \ldots \cdot \hat{\alpha}_n(\xi) = \hat{\gamma}(\xi)$. Taking logarithm gives now a linear relationship $\log \hat{\alpha}_1(\xi) + \ldots + \log \hat{\alpha}_n(\xi) = \log \hat{\gamma}(\xi)$.  Since $\alpha_k = \rho_{x_k}$,
\[
\log \hat{\rho}_{x_1}(\xi) + \ldots + \log \hat{\rho}_{x_n}(\xi) = \log \hat{\gamma}(\xi).
\]
The task is now to recover $x_1,\ldots,x_n$ from $\log \hat{\gamma}(\xi)$. This is a linear sparse inverse problem with the kernel $G(\xi,x) = \log \hat{\rho}_{x}(\xi)$ and the weights $w_k$ equal to $1$.

The second step solves this problem using the eigenmatrix method. We pick a set of samples $\xi_j\in\R$ and define for each $x\in X$ the vector
\[
\bb_x:=  \bbm \log\hat{\rho}_x(\xi_j) \ebm_j.
\]
Choose a Chebyshev grid $\{c_t\}_{1\le t \le \nc}$ over the interval $X$. If $\bb_x$ is analytic in $x$ and the columns of $\Bh$ are numerically linearly independent, we can build the matrix $M$ that satisfies $M\bb_x \approx x \bb_x$. Finally, define the vector
\[
\bu = \bbm \log\hat{\gamma}(\xi_j) \ebm_j
\]
and form the $T$ matrix to recover $x_1,\ldots, x_n$ and equivalently $\alpha_1,\ldots,\alpha_n$.

Below, we present a few numerical examples. In order to take the logarithm safely without hitting zero, we perturb the samples $\xi_j$ with an infinitesimal imaginary part.

\begin{example} The problem setup is as follows. 
  \begin{itemize}
  \item The 1-parameter family is $ \rho_x = \frac{1}{2} \delta_{0} + \frac{1}{2} \delta_x $ for $x\in [0.2, 1]$. Hence, for a fixed $x$, the distribution $\rho_x$ is supported at $0$ and $x$ with equal weights.
  \item $(\alpha_1,\alpha_2,\alpha_3)$ correspond to $(x_1,x_2,x_3) = (0.2, 0.6, 1.0)$.
  \item The number of samples from $\gamma$ is $N=102400$.
  \end{itemize}
  Here $\bb_x:= [\log\hat{\rho}_x(\xi_j)]_j $ is analytic in $x$ and the matrix $\Bh$ is verified numerically to be well-conditioned.
  Figure \ref{fig:p1} summarizes the result. The first plot is the histogram of the samples from $\gamma(y)$. The second plot shows the exact (red) and reconstructed (blue) parameters $\{x_k\}$.
  
  \begin{figure}[h!]
    \centering
    \includegraphics[scale=0.36]{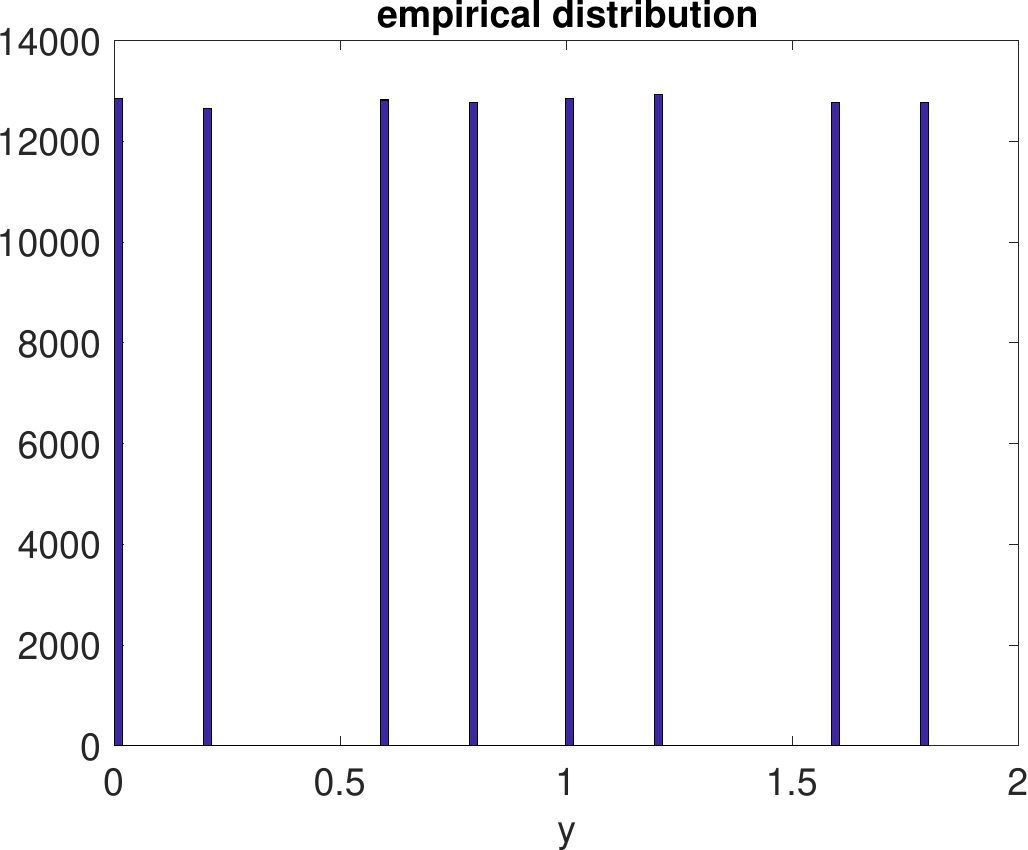}
    \includegraphics[scale=0.36]{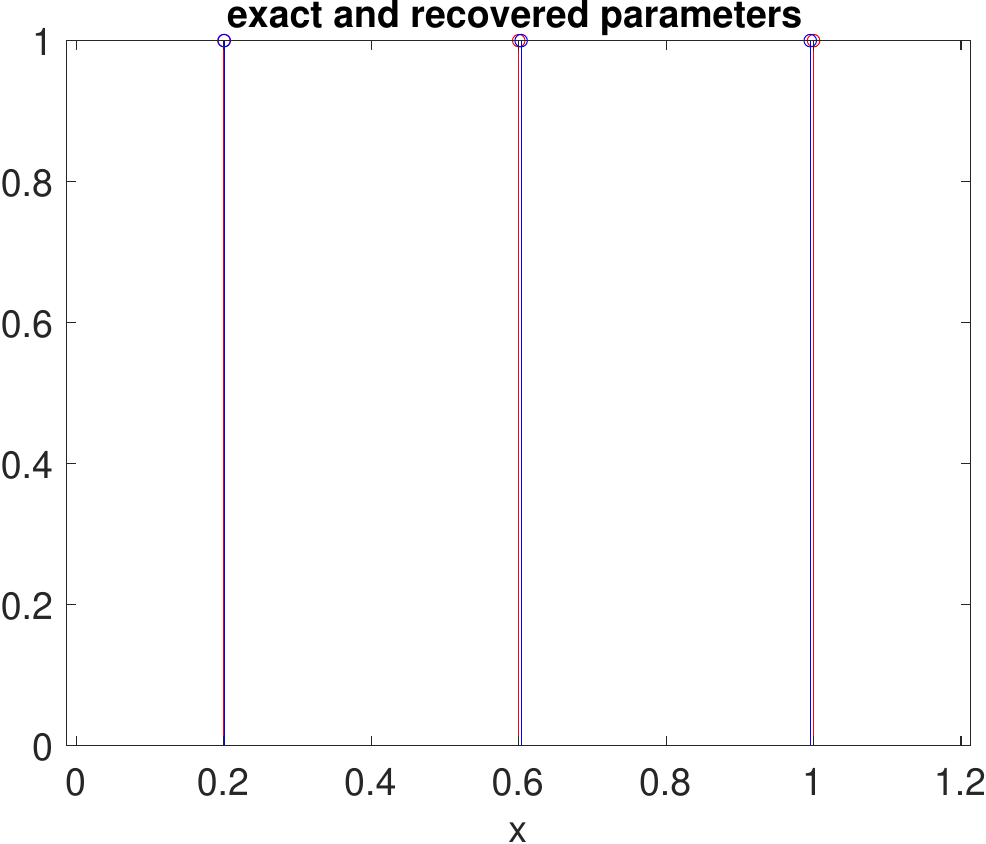}
    \caption{Left: the histogram of samples from $\gamma(y)$. Right: the exact (red) and reconstructed (blue) parameters $\{x_k\}$.}
    \label{fig:p1}
  \end{figure}

\end{example}

\begin{example} The problem setup is as follows.
  \begin{itemize}
  \item The 1-parameter family is $\rho_x = \frac{2}{3} \delta_{0} + \frac{1}{3} \delta_x$ for $x\in [0.2, 1]$.
  \item $(\alpha_1,\alpha_2,\alpha_3)$ correspond to $(x_1,x_2,x_3) = (0.2, 0.6, 1.0)$.
  \item The number of samples from $\gamma$ is $N=102400$.
  \end{itemize}
  $ \bb_x:= [\log\hat{\rho}_x(\xi_j)]_j $ is again analytic in $x$ and the matrix $\Bh$ is verified numerically to be well-conditioned.  Figure \ref{fig:p2} summarizes the result. The first plot is the histogram of the samples from $\gamma(y)$. The second plot gives the exact and reconstructed parameters $\{x_k\}$.

  \begin{figure}[h!]
    \centering
    \includegraphics[scale=0.36]{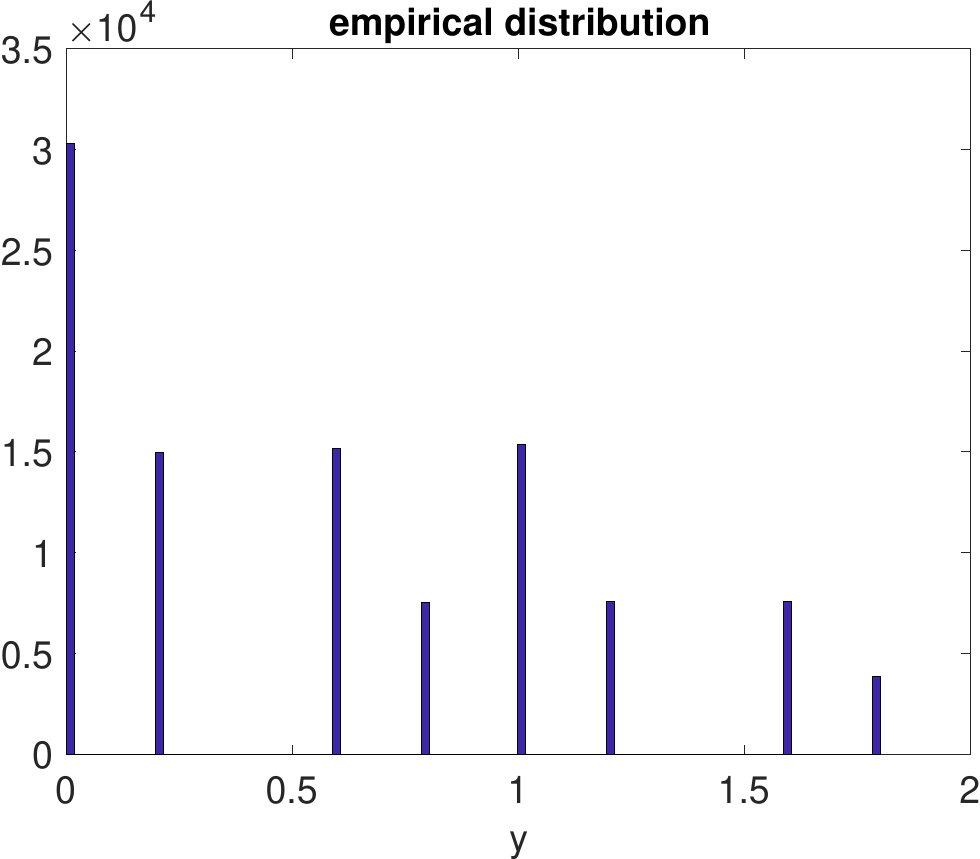}
    \includegraphics[scale=0.36]{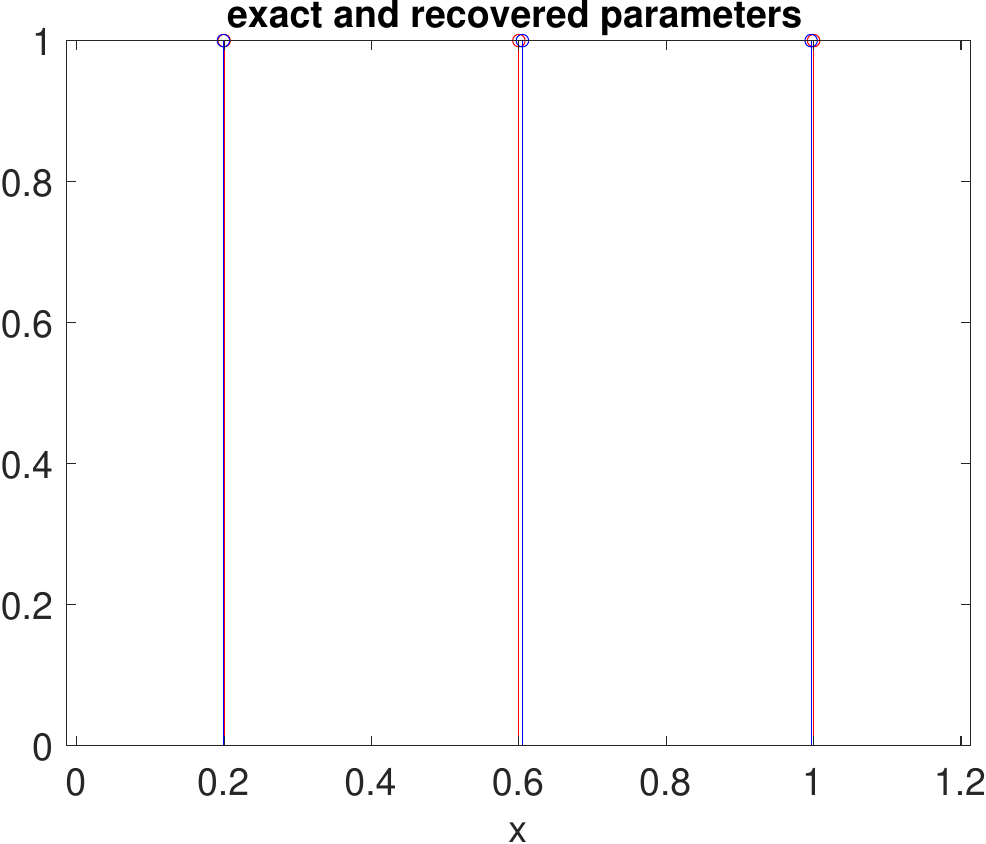}
    \caption{Left: the histogram of samples from $\gamma(y)$. Right: the exact (red) and reconstructed (blue) parameters $\{x_k\}$.}
    \label{fig:p2}
  \end{figure}

\end{example}

\section{Free additive deconvolution}\label{sec:add}

Let the interval $X$ be the parametric domain of the 1-parameter family $\rho_x$. For each $A_k$, the spectral measure $\mu_{A_k} = \rho_{x_k}$ for unknown $x_k$. Given $\muC$ of $C = A_1 + Q_2 A_2 Q_2^\T + \ldots + Q_n A_n Q_n^\T$, the task is to recover $\mu_{A_1}, \ldots, \mu_{A_n}$ or equivalently $x_1,\ldots,x_n$.

The first step turns this nonlinear free deconvolution problem into a linear one. Given any $\mu$, the Cauchy integral $g = \int \frac{1}{z-y} d\mu(y)$ defines a correspondence between $z$ and $g$. The map from $z$ to $g$ is the Stieltjes transform, denoted by $g_\mu(z)$. Its inverse is well-defined for sufficiently small values of $g$ and is denoted by $z_\mu(g)$. The R-transform \cite{voiculescu1986addition} is defined as
\[
r_\mu(g) = z_\mu(g)-\frac{1}{g}.
\]
For each $A_k$, $\mu_{A_k} = \rho_{x_k}$, and its R-transform is denoted as $r_{\rho_{x_k}}(g)$. In the large dimension limit, $r_\muC(g) = r_\muAa(g) + \ldots + r_\muAn(g)$, \cite{voiculescu1986addition} or equivalently,
\[
r_\muC(g) = r_{\rho_{x_1}}(g) + \ldots r_{\rho_{x_n}}(g).
\]
The task is, given $r_\muC(g)$, to recover $x_1,\ldots,x_n$. This is now a linear sparse inverse problem with the kernel given by $G(g,x)= r_{\rho_{x}}(g)$ and the weights $w_k$ all equal to $1$.

The second step solves this problem using the eigenmatrix method. Assume for now that for any $\mu$, we are able to compute $r_\mu(g)$. We choose a set of samples $\{g_j\}$ on a small circle around the origin. For each $x\in X$, define the vector
\[
\bb_x:=  \bbm r_{\rho_x}(g_j) \ebm_j.
\]
Choose a Chebyshev grid over the $\{c_t\}_{1\le t \le \nc}$ over the interval $X$. If $\bb_x$ is analytic in $x$ and the columns of $\Bh$ are numerically linearly independent, we can build the matrix $M$ that satisfies $M\bb_x \approx x \bb_x$. Finally, define the vector
\[
\bu = \bbm r_\muC(g_j) \ebm_j
\]
and form the $T$ matrix to recover $x_1,\ldots, x_n$, and equivalently $\mu_{A_1},\ldots,\mu_{A_n}$.

Let us address the question of, given $\mu$, how to compute its R-transform $r_\mu(g)$ for given $g$. This is used heavily in the definition of $\bb_x$ and $\bu$. It is sufficient to compute $z_\mu(g)$, which is equivalent to solving the equation $g_\mu(z) = g$. This is done with a Newton iteration. Let $z_i$ be the current approximation. In order to find the update $\eps$, we write
\[
g_\mu(z_i + \eps) \approx g.
\]
Linear expansion at $z_i$ results in $g_\mu(z_i) + g_\mu'(z_i) \eps \approx g$. Hence, we set $ \eps \equiv g_\mu'(z_i)^{-1} (g-g_\mu(z_i))$, resulting the Newton iteration
\[
z_{i+1} = z_i + g_\mu'(z_i)^{-1} (g-g_\mu(z_i)).
\]
The limit of $\{z_i\}$ gives $z_\mu(g)$ and we set $r_\mu(g) = z_\mu(g) - 1/g$.

Below, we present a few numerical examples. In each example, we select a 1-parameter family with the first moment equal to zero. The reason is that a non-zero first moment can be absorbed in any $\mu_{A_k}$ with a uniform spectral shift.

\begin{example} The problem setup is as follows.
  \begin{itemize}
  \item The 1-parameter family is $ \rho_x = \frac{2}{3} \delta_{-x/2} + \frac{1}{3} \delta_x  $ for $x\in X=[0.4, 1]$.
  \item The spectral measure of $A_1,A_2$ correspond to $(x_1,x_2) = (0.5, 0.9)$.
  \item The dimension is $N=8192$.
  \end{itemize}
  Here $ \bb_x:= [r_{\rho_x}(g_j)]_j $ is analytic in $x$. Though the condition number of $\Bh$ is quite large, the construction of $M$ is quite accurate. Figure \ref{fig:a2} summarizes the result. The first plot is the histogram of $\muC$. The second plot gives the exact (red) and reconstructed (blue) parameters $\{x_k\}$.
  \begin{figure}[h!]
    \centering
    \includegraphics[scale=0.36]{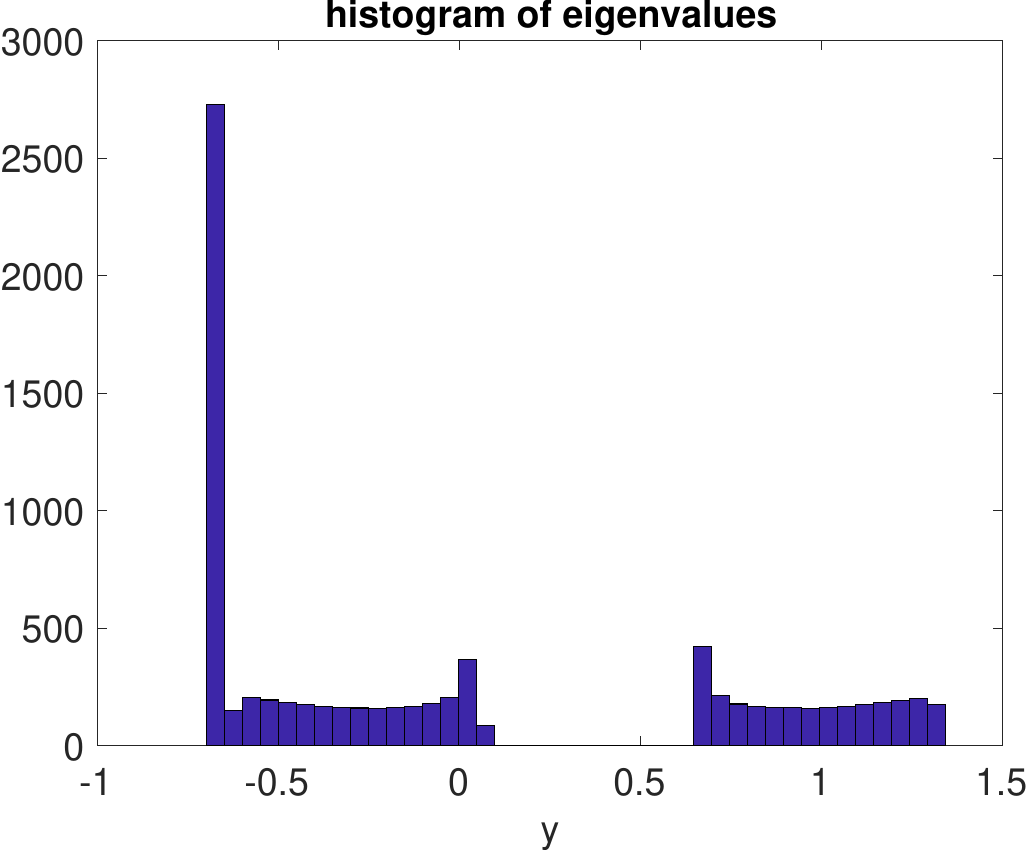}
    \includegraphics[scale=0.36]{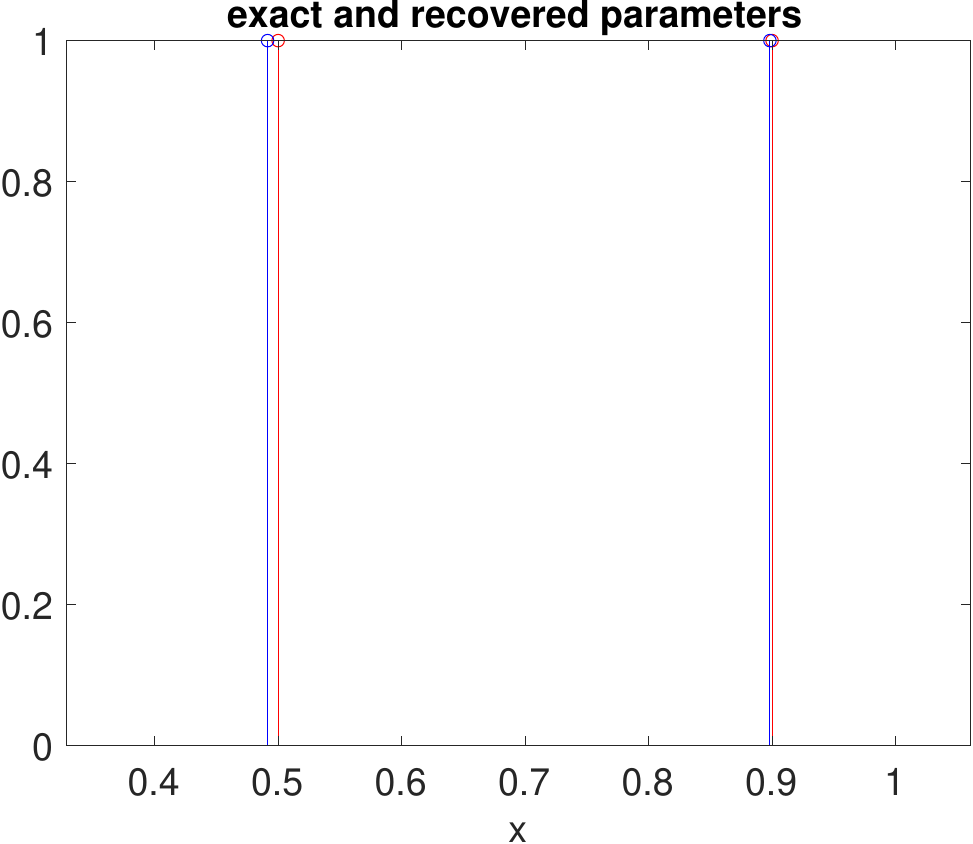}
    \caption{Left: the histogram of $\mu_C$. Right: the exact (red) and reconstructed (blue) parameters $\{x_k\}$.}
    \label{fig:a2}
  \end{figure}
\end{example}

\begin{example} The problem setup is as follows. 
  \begin{itemize}
  \item The 1-parameter family is $ \rho_x = \frac{1}{2} \delta_{-x} + \frac{1}{2} \delta_x $ for $x\in X=[0.4, 1]$.
  \item The spectral measures of $A_1,A_2,A_3$ correspond to $(x_1,x_2,x_3) = (0.4, 0.7, 1.0)$.
  \item The dimension is $N=8192$.
  \end{itemize}
  Here $ \bb_x:= [r_{\rho_x}(g_j)]_j $ is analytic in $x$. Though the condition number of $\Bh$ is quite large, the construction of $M$ is still quite close. Figure \ref{fig:a1} summarizes the result. The first plot is the histogram of $\muC$. The second plot gives the exact (red) and reconstructed (blue) parameters $\{x_k\}$.
  \begin{figure}[h!]
    \centering
    \includegraphics[scale=0.36]{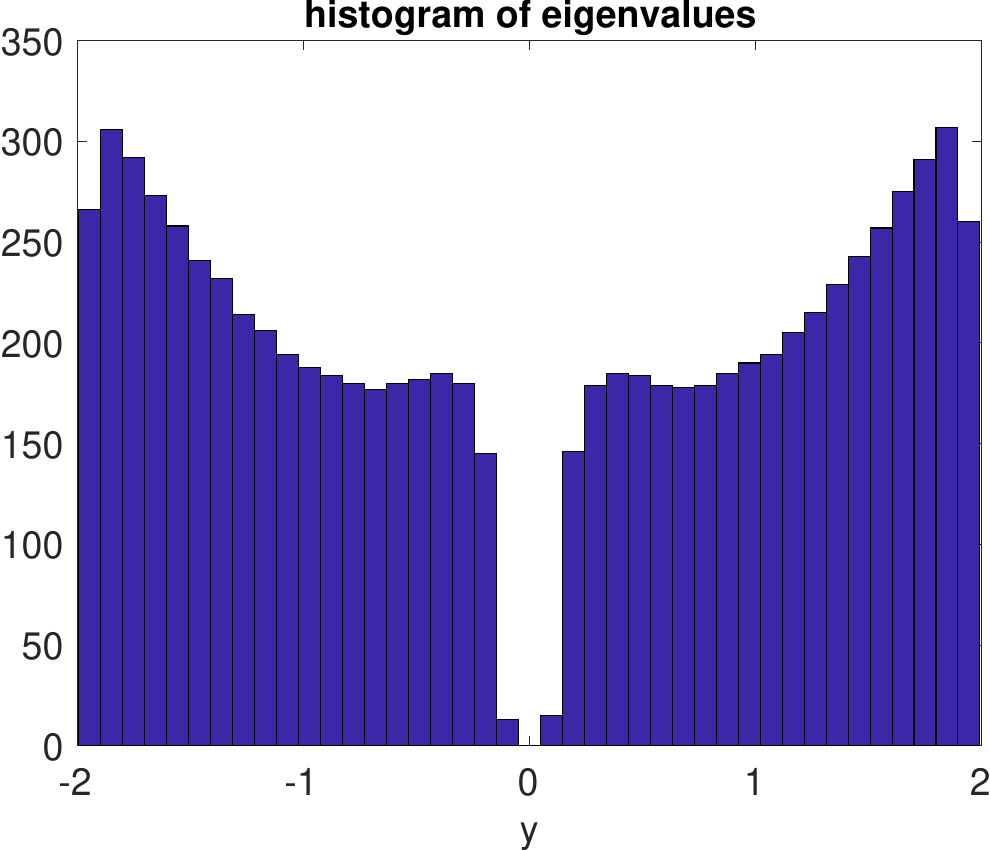}
    \includegraphics[scale=0.36]{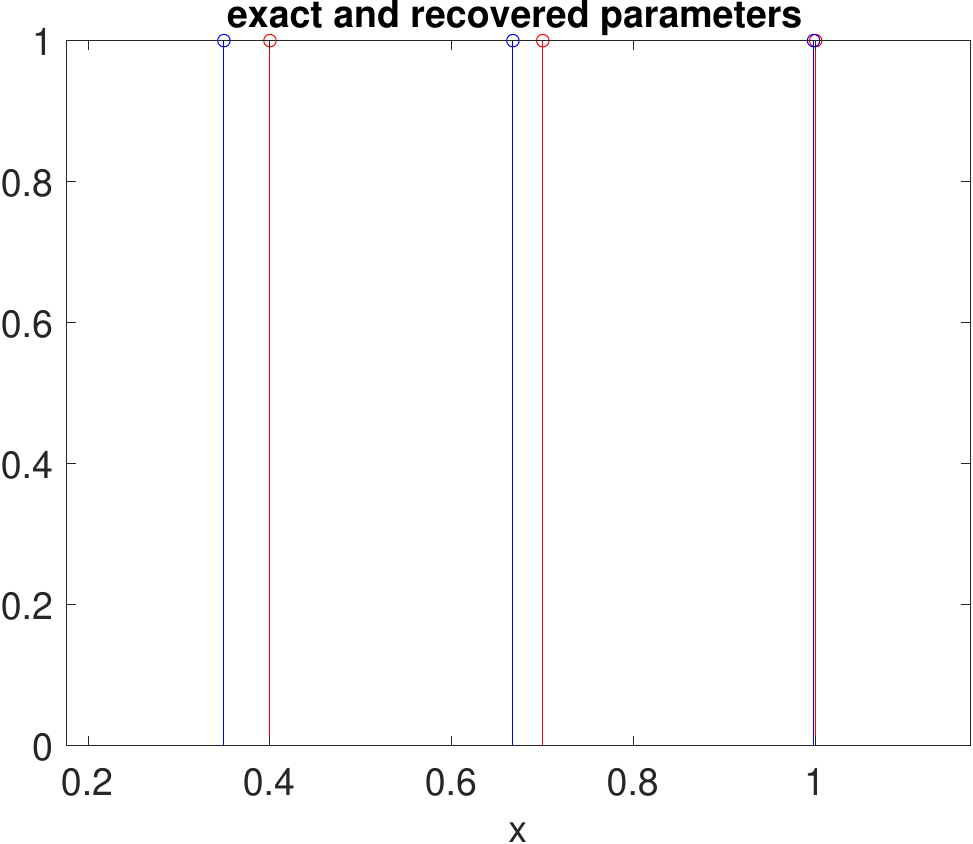}
    \caption{Left: the histogram of $\mu_C$. Right: the exact (red) and reconstructed (blue) parameters $\{x_k\}$.}
    \label{fig:a1}
  \end{figure}
\end{example}

\section{Free multiplicative deconvolution}\label{sec:mul}

Let the interval $X$ be the parametric domain of the 1-parameter family $\rho_x$. For each $k$, the spectral measure $\mu_{A_k} = \rho_{x_k}$ for unknown $x_k$.  Given the spectral measure $\muC$ of the product $C =\sqrt{A_1} Q_2 \sqrt{A_2} \ldots Q_n A_n Q_n^\T \ldots \sqrt{A_2} Q_2^\T \sqrt{A_1}$, the task is to recover $\mu_{A_1}, \ldots, \mu_{A_n}$ or equivalently $x_1,\ldots,x_n$.

The first step turns this nonlinear free deconvolution problem into a linear one. Given any $\mu$, the Cauchy integral $g = \int \frac{1}{z-y} d\mu(y)$ defines a correspondence between $z$ and $g$. In addition, introduce $t = z g-1$ and $s = \frac{t+1}{t z}$. The correspondence from $t$ to $s$ is well-defined for sufficiently small values of $t$ and is called the S-transform \cite{voiculescu1987multiplication}, denoted as $s_\mu(t)$.  For each $A_k$, $\mu_{A_k} = \rho_{x_k}$, and its S-transform is denoted as $s_{\rho_{x_k}}(t)$. In the large dimension limit, $s_\muC(t) = s_\muAa(t)\cdot \ldots \cdot s_\muAn(t)$, or equivalently
\[
s_\muC(t) = s_{\rho_{x_1}}(t)\cdot \ldots \cdot s_{\rho_{x_n}}(t).
\]
Taking logarithms leads to
\[
\log s_\muC(t) = \log s_{\rho_{x_1}}(t) + \ldots + \log s_{\rho_{x_n}}(t).
\]
The task is, given $\log s_\muC(t)$, to recover $x_1,\ldots,x_n$. This is a linear sparse inverse problem with the kernel $G(t,x)= \log s_{\rho_{x}}(t)$ and all weights $w_k$ equal to 1.

The second step solves this problem using the eigenmatrix method. Assume for now that for any $\mu$, we are able to compute $s_\mu(t)$.  Choose a set of samples $\{t_j\}$ on a small circle around the origin. For each $x\in X$, define the vector
\[
\bb_x:=  \bbm \log s_{\rho_x}(t_j) \ebm_j.
\]
Choose a Chebyshev grid over the $\{c_t\}_{1\le t \le \nc}$ over the interval $X$. If $\bb_x$ is analytic in $x$ and the columns of $\Bh$ are numerically linearly independent, we can build the matrix $M$ that satisfies $M\bb_x \approx x \bb_x.$ Finally, define the vector
\[
\bu = \bbm \log s_\muC(t_j) \ebm_j
\]
and form the $T$ matrix to recover $x_1,\ldots,x_n$, and equivalently $\mu_{A_1},\ldots,\mu_{A_n}$.

Let us now address the question of given $\mu$, how to compute its S-transform $s_\mu(t)$ for given $t$. This is used in the definition of $\bb_x$ and $\bu$.
From the definition of the S-transform, the correspondence among $z, g, t, s$ is given by
\[
t+1 = z g, \quad s = (t+1)/(tz).
\]
Given $t$, one way to compute $s$ is to recover $z$ first, via solving the equation $ t+1 = z g_\mu(z)$.  This is again done with the help of a Newton iteration. Let $z_i$ be the current approximation. To find the update $\eps$, we write
\[
(z_i+\eps) g_\mu(z_i + \eps) = t+1.
\]
Linear expansion at $z_i$ gives $z_ig_\mu(z_i) + [g_\mu(z_i) + z_i g_\mu'(z_i)] \eps \approx t+1$. Hence, we set $\eps = (g_\mu(z_i) + z_i g_\mu'(z_i))^{-1} (t+1 - z_i g_\mu(z_i)).$ This gives the Newton iteration
\[
z_{i+1} = z_i + (g_\mu(z_i) + z_i g_\mu'(z_i))^{-1} (t+1 - z_i g_\mu(z_i)).
\]
The limit of $\{z_i\}$ gives $z_\mu(t)$. Then $s_\mu(t) = (t+1)/(t z_\mu(t))$ follows easily.

Below, we present a few numerical examples. In each example, the 1-parameter family has the first moment equal to one. This is because a non-identity uniform spectral scaling can be absorbed in either $\mu_{A_1}$ or $\mu_{A_2}$.

\begin{example} The problem setup is as follows. 
  \begin{itemize}
  \item The 1-parameter family is $ \rho_x = \frac{2}{3} \delta_{3/(2+x)} + \frac{1}{3} \delta_{3x/(2+x)}$ for $x\in X=[1.4, 3]$.
  \item The spectral measures of $A_1,A_2$ correspond to $(x_1,x_2) = (1.7, 2.5)$.
  \item The dimension is $N=8192$.
  \end{itemize}
  $ \bb_x:= [\log s_{\rho_x}(t_j)]_j $ again is analytic in $x$. Though the condition number of $\Bh$ is quite large, the construction of $M$ is quite accurate.  Figure \ref{fig:m2} summarizes the result. The first plot is the histogram of $\muC$. The second plot gives the exact (red) and reconstructed (blue) parameters $\{x_k\}$.
  \begin{figure}[h!]
    \centering
    \includegraphics[scale=0.36]{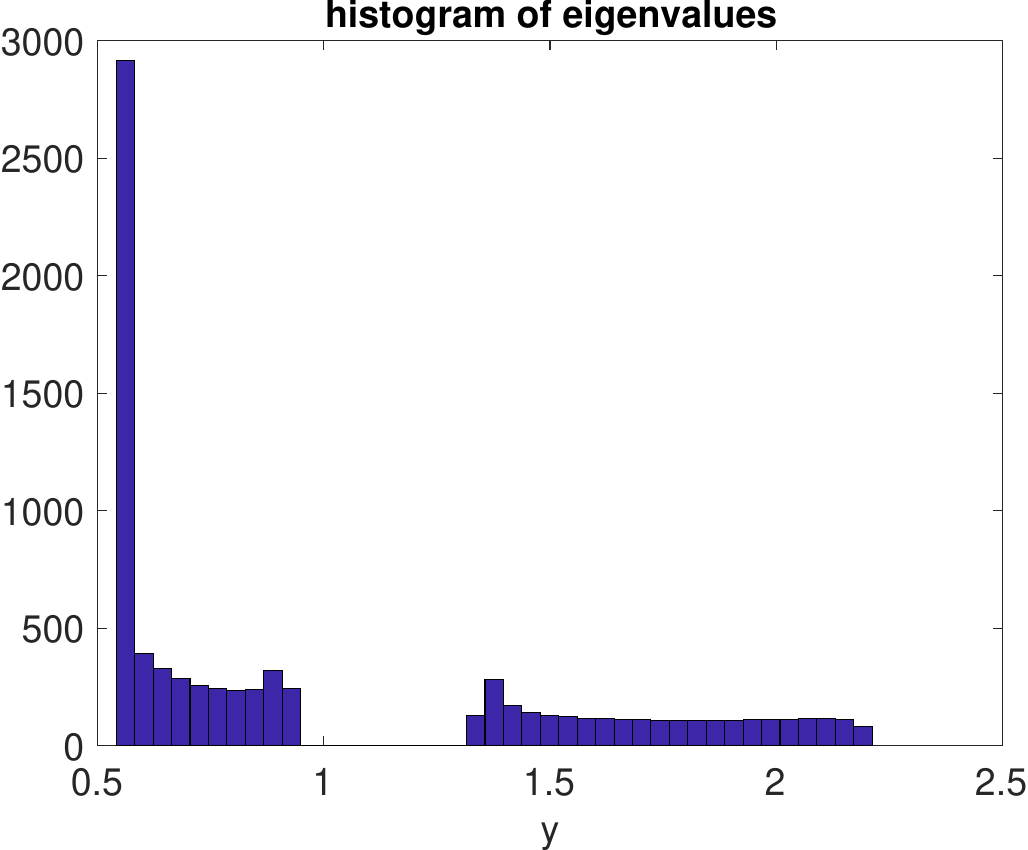}
    \includegraphics[scale=0.36]{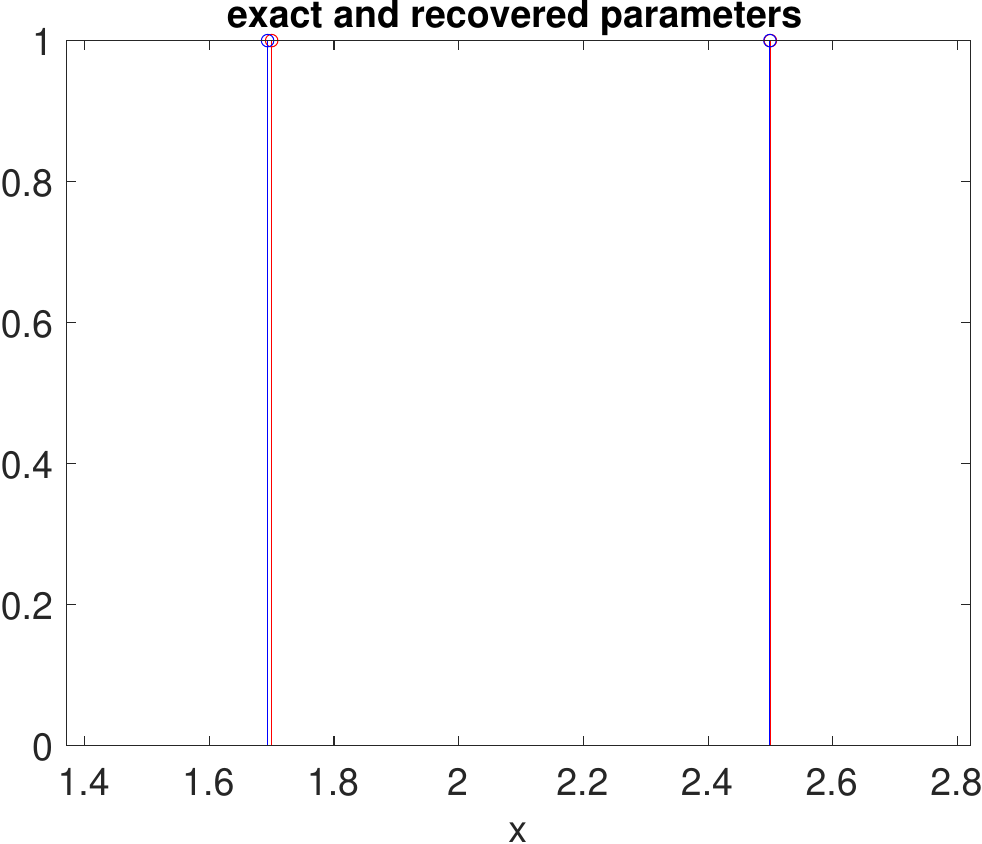}
    \caption{Left: the histogram of $\mu_C$. Right: the exact (red) and reconstructed (blue) parameters $\{x_k\}$.}
    \label{fig:m2}
  \end{figure}  
\end{example}

\begin{example} The problem setup is as follows. 
  \begin{itemize}
  \item The 1-parameter family is $ \rho_x = \frac{1}{2} \delta_{2/(1+x)} + \frac{1}{2} \delta_{2x/(1+x)}$ for $x\in X=[1.4, 3]$.
  \item The spectral measures of $A_1,A_2,A_3$ correspond to $(x_1,x_2,x_3) = (1.4, 2.2, 3.0)$.
  \item The dimension is $N=8192$.
  \end{itemize}
  Here $ \bb_x:= [\log s_{\rho_x}(t_j)]_j $ is analytic in $x$. Though the condition number of $\Bh$ is quite large, the construction of $M$ is still quite close.  Figure \ref{fig:m1} summarizes the result. The first plot is the histogram of $\muC$. The second plot gives the exact (red) and reconstructed (blue) parameters $\{x_k\}$.
  \begin{figure}[h!]
    \centering
    \includegraphics[scale=0.36]{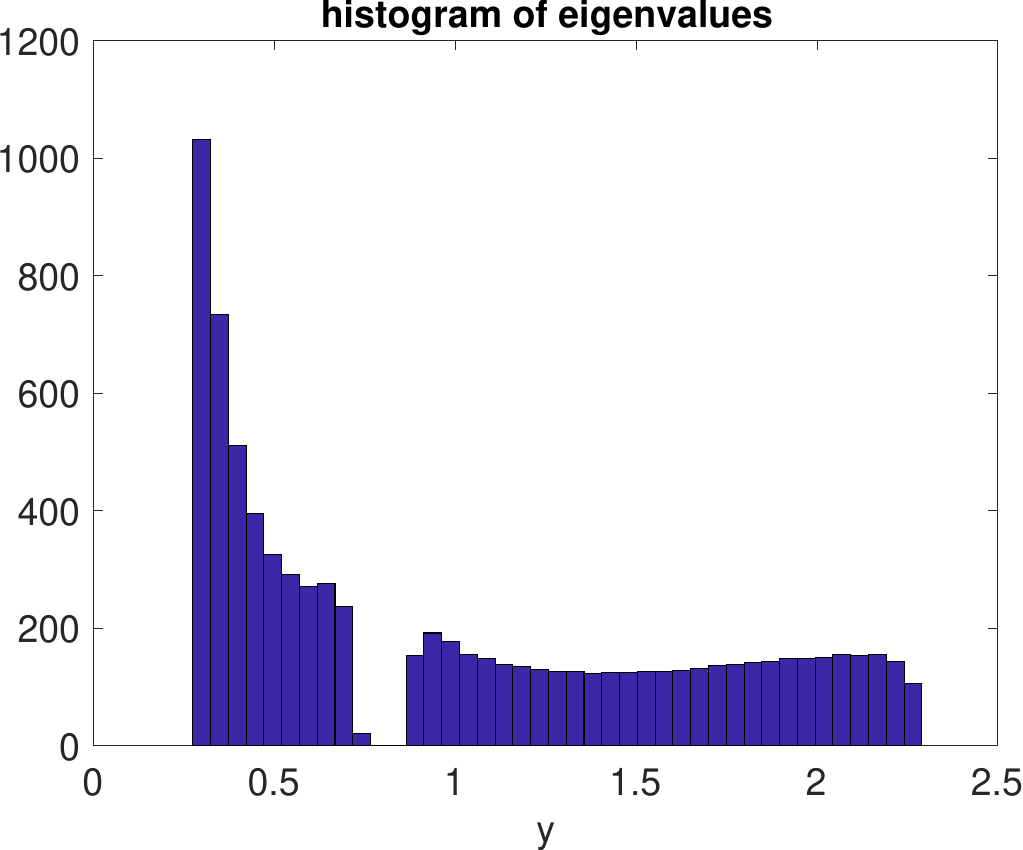}
    \includegraphics[scale=0.36]{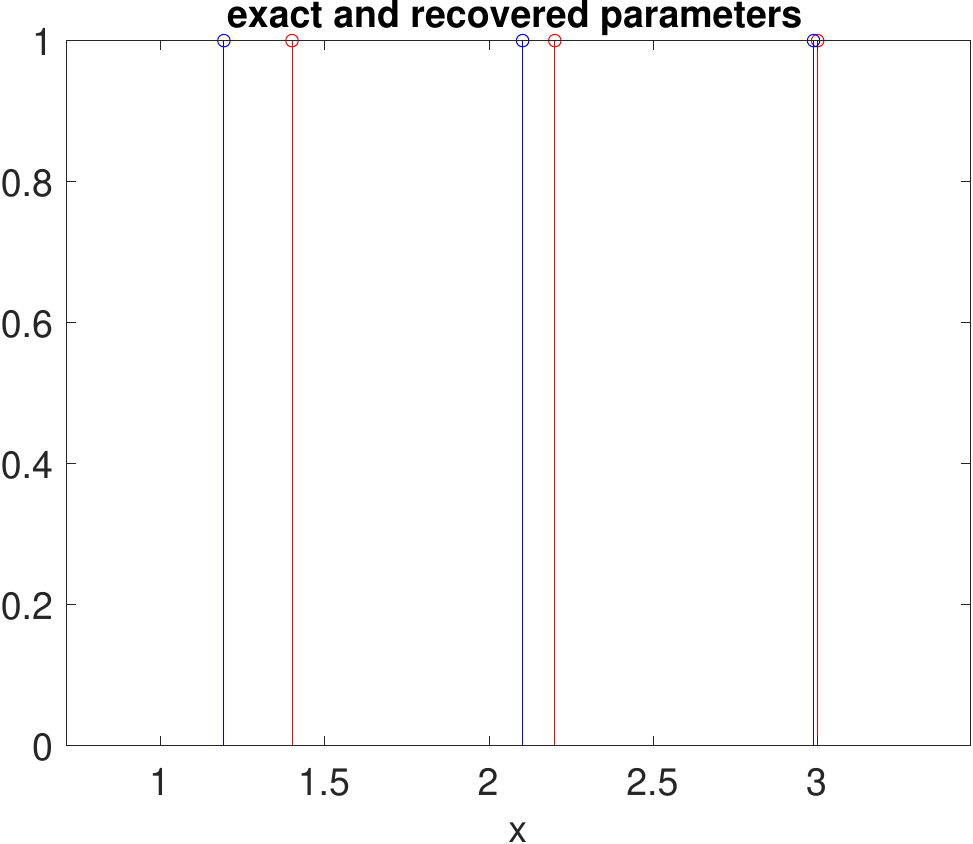}
    \caption{Left: the histogram of $\mu_C$. Right: the exact (red) and reconstructed (blue) parameters $\{x_k\}$.}
    \label{fig:m1}
  \end{figure}
\end{example}

\section{Discussion}\label{sec:disc}

This note is a preliminary exploration of blind free deconvolution problems. Below are several directions for improvements and future work.

First, we only consider 1-parameter families for the spectral measure of $A_k$. This is clearly restrictive, and an obvious direction is to include families with multiple parameters. The work of \cite{ying2025multidimensional} on multidimensional structure sparse recovery can potentially be leveraged.

Second, the sparsity assumption of the spectral measure of $A$ may not be appropriate in many cases. For a specific application, if the spectral measure of $A$ arises from some other low-complexity models controlled by a small number of parameters, this approach based on an eigenmatrix should also apply.

Third, the eigenmatrix method can be viewed as a non-Hermitian quantization of the spectral variety $X$. The vectors $\bb_x$ serve as the approximate non-orthogonal eigenvectors. Once the matching condition $M \bb_x \approx x \bb_x$ is satisfied on a sufficiently dense Chebyshev grid, the whole spectral variety $X$ resides in the pseudo-spectrum of $M$. The original eigenmatrix method was proposed for simple $X$ (either an interval or a domain in the complex plane). The work here is a first step towards addressing a more general spectral variety, even though we consider only a curve within it by focusing on a one-parameter family.

\bibliographystyle{abbrv}

\bibliography{ref}

\end{document}